\documentclass[10pt]{amsart}
\usepackage{amsmath, amsfonts, amssymb, amstext, amscd, amsthm}

\newtheorem{algebra1}{Theorem}[section]
\newtheorem{algebra2}[algebra1]{Theorem}
\newtheorem{algebra3}[algebra1]{Theorem}

\newtheorem{ternary1}[algebra1]{Definition}
\newtheorem{ternary2}[algebra1]{Definition}

\newtheorem{fractions1}[algebra1]{Proposition}
\newtheorem{fractions2}[algebra1]{Proposition}

\newtheorem{cor1}[algebra1]{Proposition}
\newtheorem{cor2}[algebra1]{Corollary}

\begin{document}


\title{A Ternary Algebra \\ with Applications to Binary Quadratic Forms}

\author{Edray Goins}

\address{Mathematical Sciences Research Institute \\ 1000 Centennial
Drive \# 5070 \\ Berkeley, CA 94720-5070}

\email{goins@msri.org}


\begin{abstract}
We discuss multiplicative properties of the binary quadratic form $a
\, x^2 + b \, x \, y + c \, y^2$ by considering a ring of matrices
which is closed under a triple product.  We prove that the ring forms
a ternary algebra in the sense of Hestenes, and then derive both
multiplicative formulas for a large class of binary quadratic forms
and a type of multiplication for points on a conic section which
generalizes the algebra of rational points on the unit circle.
\end{abstract}

\footnote{This research was sustained in part by fellowship stipend support
from the National Physical Science Consortium and the National
Security Agency.} 

\maketitle


\section{Introduction}

The multiplicative relation
\[ \left[ x_1^2 + y_1^2 \right] \left[ x_2^2 + y_2^2 \right]  = \left(
x_1 \, x_2 + y_1 \, y_2 \right)^2 + \left( y_1 \, x_2 - x_1 \, y_2 \right)^2
\]

\noindent is well-known.  This and similar formulas were known to Euler;
see \cite{MR90m:11016}.  This exposition is motivated by the triple product
\[ \left[ a \, x_1^2 + b \, x_1 \, y_1 + c \, y_1^2 \right] \left[ a \,
x_2^2 + b \, x_2 \, y_2 + c \, y_2^2 \right] \left[ a \, x_3^2 + b \,
x_3 \, y_3 + c \, y_3^2 \right] \]

\noindent which is also in the form $a \, x^2 + b \, x \, y + c \, y^2$
where
\[ \begin{aligned}
   x & = a \, (x_1 \, x_2 \, x_3) + b \, (x_1 \, y_2 \, x_3) + c \,
   (x_1 \, y_2 \, y_3 - y_1 \, x_2 \, y_3 + y_1 \, y_2 \, x_3) \\
   y & = a \, (x_1 \, x_2 \, y_3 - x_1 \, y_2 \, x_3 + y_1 \, x_2 \,
   x_3) + b \, (y_1 \, x_2 \, y_3) + c \, (y_1 \, y_2 \, y_3)
\end{aligned} \]

I consider the multiplicative properties of the most general binary
quadratic form $q(x,y) = a \, x^2 + b \, x \, y + c \, y^2 + d \, x + e
\, y + f$ as defined over $\mathbb Z$.  The results are divided into three
sections.  (The results contained within this paper are stated for $\mathbb
Q$, but they hold for any field of characteristic different from 2.)

First, I consider finite-dimensional linear spaces $\mathcal A$ which
are closed under a triple product $A \cdot B^* \cdot C$ for $A,B,C \in
\mathcal A$.  
Such a class is called a ternary algebra, and was first studied in
\cite{MR25:1172}.  (Consider \cite{MR25:456}, \cite{MR27:169} and
\cite{MR56:3174} as well.)  I show that each such linear space is intimately
connected with a commutative ring $\mathcal R$ which has a nontrivial left
action on $\mathcal A$.

Second, I focus on conic sections $\mathcal C(q)$ defined by
$q(x,y)=0$ and show how they are associated with a specific ternary algebra
$\mathcal A(q)$.  I use the multiplicative structure of the elements
in $\mathcal A(q)$ with a specific fixed determinant to define a
multiplicative structure on $\mathcal C(q)$, which generalizes the
multiplicative structure of rational points on the unit circle.  This
motivates 
the concept of a ternary group, and I show the existence of an abelian group
which acts transitively.

Third, I show how the determinant map $\mathcal A(q) \to \mathbb Q$ yields
the multiplicative 
formulas stated above. I conclude by considering the rational values of
the quadratic form $q(x,y)$ and showing that they form a ternary group as
well.


\section{Ternary Algebras}

Throughout we fix $q(x,y) = a \, x^2 + b \, x \, y + c \, y^2 + d \, x + e
\, y + 
f \in \mathbb Z[x, \, y]$ and define
\begin{equation} \label{constants}
   \text{Det}(q) = \left| \begin{matrix} a & b/2 & d/2 \\ b/2 & c & e/2 \\
   d/2 & e/2 & f \end{matrix} \right|, \qquad \text{Disc}(q) = \left|
   \begin{matrix} a & b/2 \\ b/2 & c \end{matrix} \right|;
\end{equation}

\noindent When $\text{Disc}(q) \neq 0$ we also define
\begin{equation} \label{center}
   h_q = \frac {\left| \begin{matrix} b/2 & d/2 \\ c & e/2 \end{matrix}
   \right|}{\text{Disc}(q)}, \qquad k_q = - \frac {\left| \begin{matrix}
   a & d/2 \\ b/2 & e/2 \end{matrix} \right|}{\text{Disc}(q)}, \qquad
   \text{and} \qquad m_q = -\frac {\text{Det}(q)}{\text{Disc}(q)}.
\end{equation}

\noindent We may identify the point $(h_q, k_q)$ with the center of
the conic section $q(x,y) = 0$.  Indeed, it is easy to verify that
\begin{equation} \label{m}
   q(x, \, y) = a \, (x - h_q)^2 + b \, (x-h_q) \, (y-k_q) + c \,
   (y-k_q)^2 - m_q
\end{equation}

\begin{algebra1} \label{algebra1} Fix $q(x,y)$ with $\text{Disc}(q) \neq
0$.  With notation as above, define
\begin{equation} \begin{aligned} \label{groups}
   \mathcal A(q) & = \left \{ \left. \left[x, \, y \right]_q = \left[
   \begin{matrix} (x-h_q) & - c \, (y-k_q) \\ (y-k_q) & a \, (x-h_q) +
   b \, (y-k_q) \end{matrix} \right] \right| \ x, \, y \in \mathbb
   Q \right \} \\
   \mathcal R(q) & = \left \{ \left. \left[u, \, v \right]_{q,0} = \left[
   \begin{matrix} u & - c \, v \\ a \, v & u + b \,
   v \end{matrix} \right] \right| \ u, \, v \in \mathbb Q \right \} \\
\end{aligned} \end{equation}

Then we have the following.
\begin{enumerate}
   \item $\mathcal A(q)$ is a linear vector space.
   \item $\mathcal R(q)$ is a commutative ring with identity.
   \item Given $A,B \in \mathcal A(q)$, the product $A \cdot B^* \in
   \mathcal R$ where $B^*$ is the transpose of the cofactor matrix of $B$.
   \item $\mathcal A(q)$ is a left $\mathcal R(q)$-module.
   \item There exist $B_0, C_0 \in \mathcal A(q)$ such that $A \cdot
   B_0^* \cdot C_0 = A$ for all $A \in \mathcal A(q)$.
\end{enumerate}

Moreover, for $\alpha \in \mathbb Q$ and $A,B,C,D,E \in \mathcal A(q)$,
\begin{enumerate}
   \item Multiplication: $A \cdot B^* \cdot C \in \mathcal A(q)$
   \item Commutativity: $A \cdot B^* \cdot C = C \cdot B^* \cdot A$
   \item Associativity: $\left( A \cdot B^* \cdot C \right) \cdot
   D^* \cdot E = A \cdot B^* \cdot \left( C \cdot D^* \cdot
   E \right)$
   \item Distributivity: $(A+D)\cdot B^* \cdot C = A\cdot B^*
   \cdot C + D \cdot B^* \cdot C$
   \item Linearity: $(\alpha \, A) \cdot B^* \cdot C = \alpha \, (A
   \cdot B^* \cdot C) = A \cdot B^* \cdot (\alpha \, C)$
   \item Nondegeneracy: $A \cdot A^* \cdot A = 0$ if and only if
   $\det A=0$
\end{enumerate}
\end{algebra1}

The motivation behind studying these matrices comes from the
determinant map: by \eqref{m}, $\det [x,y]_q= q(x, \, y) + m_q$.  Once we
know the multiplicative properties of $\mathcal A(q)$ then we may deduce
similar properties for $q(x,y)$.  Note that in general $\mathcal A(q)$
is not closed under multiplication so the triple product is indeed
necessary.  We remind the reader that the transpose of the cofactor matrix
satisfies $B^* \cdot B = \det B \cdot 1$, and is explicitly given by
\[ \left[x, \, y \right]_q^* = \left[ \begin{matrix} a \, (x-h_q) + b
\, (y-k_q) & c \, (y-k_q) \\ - (y-k_q) & (x-h_q) \end{matrix} \right] \]

\begin{proof}  It is obvious that $\mathcal A(q)$ is a linear vector
space.  To prove that $\mathcal R(q)$ is a ring, we note that $\mathcal
R(q)$ is a subset of the ring $\text{Mat}_2(\mathbb Q)$ so it suffices to
show that for all $g,h \in \mathcal R(q)$ we have $g-h$, $g \cdot h$, and
the multiplicative  identity all in $\mathcal R$.  Clearly $[1, \,
0]_{q,0} \in \mathcal R(q)$ is the multiplicative identity, and, denoting
$g=[u_1,v_1]_{q,0}$ and $h=[u_2,v_2]_{q,0}$, we see that the difference $g-h
= 
\left[u_1-u_2, \, v_1-v_2 \right]_{q,0}$ and the product $g \cdot h = \left[
u_1 \, u_2 - a \, c \, v_1 \, v_2, \ \ u_1 \, v_2 + v_1 \, u_2 + b \,
u_1 \, v_2 \right]_{q,0}$ are indeed in $\mathcal R(q)$.  $g \cdot h =
h\cdot 
g$ is clear from this expression so $\mathcal R(q)$ is commutative.

One calculates explicitly that $[u,v]_{q,0} = [x_1,y_1]_q \cdot
[x_2,y_2]_q^* \in \mathcal R(q)$ where
\[ \begin{aligned}
   u & = a \, (x_1-h_q) \, (x_2-h_q) + b \, (x_1-h_q) \, (y_2-k_q) + c
   \, (y_1-k_q) \, (y_2-k_q) \\
   v & = (x_2-h_q) \, (y_1-k_q) - (x_1-h_q) \, (y_2-k_q)
\end{aligned} \]

\noindent Conversely, given $g=[u,v]_{q,0} \in \mathcal R(q)$, $h \in
\mathcal R(q)$ and $C=[x,y]_q
\in \mathcal A(q)$ one also verifies that $g \cdot C = \left[ u \, x -
c \, v \, y, \ \ u \, y + a \, v \, x + b \, v \, y \right]_q$ is indeed
in $\mathcal A(q)$.  The expressions $(g \cdot h) \cdot C = g \cdot (h
\cdot C)$ and $(g + h) \cdot C = g \cdot C + h \cdot C$ follow from
the associativity and distributivity of $\text{Mat}_2(\mathbb Q)$ proving
that $\mathcal A(q)$ is indeed a left $\mathcal R(q)$-module.
   
To show the existence of elements $B_0, C_0 \in \mathcal A(q)^{\times}$ such
that $A \cdot B_0^* \cdot C_0 = A$, we note that $\text{Disc}(f)$ is nonzero
so there are rational $x_0,y_0$ such that $\det [x_0,y_0]_q \neq 0$.  Now
choose $B_0=[x_0,y_0]_q$ and $C_0 = \frac 1{\det B_0} B_0$ so that $B_0^*
\cdot 
C_0$ is the identity matrix.  The identity $A \cdot B_0^* \cdot C_0 =
A$ follows.

To show that $A\cdot B^* \cdot C \in \mathcal A(q)$ note that $g = A
\cdot B^* \in \mathcal R(q)$, and so $g \cdot C \in \mathcal A(q)$ since
$\mathcal A(q)$ is a left $\mathcal R(q)$-module.  Commutativity
comes from that of $\mathcal R(q)$:
\[ \begin{aligned}
   A \cdot B^* \cdot C & = \left( A \cdot B_0^* \cdot C_0 \right)
   \left(B \cdot B_0^* \cdot C_0 \right)^* \left( C \cdot B_0^* \cdot
   C_0 \right) \\ & = \left( A \cdot B_0^* \right) \left(C_0 \cdot C_0^*
   \right) \left( B_0 \cdot B^* \right) \left( C \cdot B_0^* \right)
   C_0 \\ & = \left( C \cdot B_0^* \right) \left(C_0 \cdot C_0^*
   \right) \left( B_0 \cdot B^* \right) \left( A \cdot B_0^* \right)
   C_0 = C \cdot B^* \cdot A
\end{aligned} \]

Associativity, distributivity, and linearity are all obvious.
Nondegeneracy follows from the relation $A \cdot A^* \cdot A = (\det
A) \, A$.  \end{proof}

The vector space $\mathcal A(q)$ is part of a larger class of spaces
which are characterized by the triple product.

\begin{ternary1} A finite-dimensional vector space $\mathcal A$ over
$\mathbb Q$ is 
called a ternary algebra if for $\alpha \in \mathbb Q$ and $A,B,C,D,E \in
\mathcal A$, there is a triple product $A \cdot B^* \cdot C$ such that
\begin{enumerate}
   \item Multiplication: $A \cdot B^* \cdot C \in \mathcal A$
   \item Commutativity: $A \cdot B^* \cdot C = C \cdot B^* \cdot A$
   \item Associativity: $\left( A \cdot B^* \cdot C \right) \cdot
   D^* \cdot E = A \cdot B^* \cdot \left( C \cdot D^* \cdot
   E \right)$
   \item Distributivity: $(A+D)\cdot B^* \cdot C = A\cdot B^*
   \cdot C + D \cdot B^* \cdot C$
   \item Linearity: $(\alpha \, A) \cdot B^* \cdot C = \alpha \, (A
   \cdot B^* \cdot C) = A \cdot B^* \cdot (\alpha \, C)$
   \item Nondegeneracy: $A \cdot A^* \cdot A = 0$ if and only if
   $\det A=0$
\end{enumerate}
\end{ternary1}

Such spaces were considered in detail in \cite{MR25:1172}, \cite{MR25:456},
\cite{MR27:169} and \cite{MR56:3174}.  \ref{algebra1} shows that the
matrices $[x,y]_q$ form
a ternary algebra.  In the proof we constructed a ring $\mathcal R(q)$
from which all of the properties of $\mathcal A(q)$ can be deduced.
Conversely, given a ternary algebra there are times when such a ring
always exists.

\begin{fractions1} \label{fractions1} Let $\mathcal A$ be a ternary algebra,
and assume 
that there exist $B_0, C_0 \in \mathcal A^{\times}$ such that $A \cdot
B_0^* \cdot C_0 = A$ for all $A \in \mathcal A$.  Denote $\mathcal R$ as
the collection of formal symbols $A \cdot B^*$ for $A,B \in \mathcal A$
modulo the equivalence relation $A \cdot B^* \sim C \cdot D^*$
whenever $A \cdot B^* \cdot C_0 = C \cdot D^* \cdot C_0$ in $\mathcal A$.

Then $\mathcal R$ is a commutative ring with identity, and $\mathcal
A$ is a left $\mathcal R$-module. \end{fractions1}

The elements $B_0$ and $C_0$ can be considered the analogue of a
multiplicative identity for ternary algebras.  Indeed, rings need not have
a multiplicative identity, but many results from commutative algebra impose
this condition.  The reader may note the similarity of the
construction of $\mathcal R$ with the construction of a field of
fractions.

\begin{proof} To show that $\mathcal R$ is a ring, we must define
addition and multiplication.  To this end, denote $\phi(B) = C_0 \cdot
B^* \cdot C_0$.  Then by the associativity of $\mathcal A$
\[ A \cdot B^* \cdot C_0 = \left(A \cdot B_0^* \cdot C_0 \right) B^* \left(
C_0 \cdot B_0^* \cdot C_0 \right) = A \cdot B_0^* \cdot \phi(B) \cdot
B_0^* \cdot C_0 \]

\noindent which shows that $A \cdot B^* \sim \left( A \cdot B_0^* \cdot
\phi(B) \right) \cdot B_0^*$.  Hence $\mathcal R$ is generated by symbols
in the form $A \cdot B_0^*$ for $A \in \mathcal A$.

Let $g = A \cdot B_0^*$ and $h = B \cdot B_0^*$ be elements of
$\mathcal R$.  Define addition by $g + h = (A + B) \cdot B_0^*$
and multiplication by $g \cdot h = (A \cdot B_0^* \cdot B) \cdot B_0^*$.
Multiplication is commutative because $g \cdot h = (A \cdot B_0^* \cdot
B) \cdot B_0^* = (B \cdot B_0^* \cdot A) \cdot B_0^* = h \cdot g$; and
it is easy to check that $1 = C_0 \cdot B_0^*$ is the multiplicative
identity.  
Associativity and distributivity follow directly from the properties of
the ternary algebra $\mathcal A$.  Hence $\mathcal R$ is a commutative
ring with identity.

Given $g = A \cdot B_0^* \in \mathcal R$ and $C \in \mathcal A$ the
definition $g \cdot C = A \cdot B_0^* \cdot C \in \mathcal A$ shows
that $\mathcal R$ acts on $\mathcal A$.  Since $\mathcal A$ is a
linear space and $\mathcal R$ is a ring, $\mathcal A$ is a left $\mathcal
R$-module.
\end{proof} 

One can use this construction to show that the symbols $B^* \cdot C$
also form a commutative ring with identity canonically
isomorphic to $\mathcal R$ which acts on $\mathcal A$ on the right.


\section{Ternary Groups}

We now use the properties of $\mathcal A(q)$ to study the points on a
conic section.

\begin{algebra2} \label{algebra2} Let $\mathcal C(q)$ be curve
defined by $q(x,y)=0$, where $\text{Det}(q), \text{Disc}(q)$ are nonzero
and $\mathcal C(q)$ has at least one rational point.  Then for $P,Q,R,S,T
\in \mathcal C(q)$, there
is a map denoted by $P \cdot Q^* \cdot R$ such that

\begin{enumerate}
   \item Multiplication: $P \cdot Q^* \cdot R \in \mathcal C(q)$
   \item Commutativity: $P \cdot Q^* \cdot R = R \cdot Q^* \cdot P$
   \item Associativity: $\left( P \cdot Q^* \cdot R \right) \cdot
   S^* \cdot T = P \cdot Q^* \cdot \left( R \cdot S^* \cdot T \right)$
   \item Identity: There exist $Q_0, R_0 \in \mathcal C(q)$ such that $P
   \cdot Q_0^* \cdot R_0 = P$ for all $P$.
   \item Inverses: Given $P$, there is $R$ such that $P \cdot Q_0^* \cdot R
= R_0$. 
\end{enumerate}
\end{algebra2}

The reader may note the similarity here with elliptic curves over $\mathbb
Q$; consult \cite{MR87g:11070}.  In order to define the group law, we must
know that a curve defined by a nondegenerate cubic does indeed contain a
rational point.

\begin{proof} First we explicitly construct this map.  Denote $\mathcal
A(q)$ as in \eqref{groups} and canonically map $\mathcal C(q) \to \mathcal
A(q)$ by $(x,y) \mapsto [x,y]_q$.  By \eqref{m} we have $\det
[x,y]_q = q(x,y) + m_q$ nonzero by assumption, so we actually have a
one-to-one map of $\mathcal C(q)$ with those matrices in $\mathcal A(q)$
with determinant $m_q$.  By the same proof of \ref{algebra1}, we see that
the product $[x_1,y_1]_q \cdot [x_2,y_2]_q^{-1} \cdot [x_3,y_3]_q = [x,y]_q$
is 
also an element of $\mathcal A(q)$.  Define $P \cdot Q^* \cdot R$ as the
point $(x,y)$.  
We must show that it is also a point on the conic section, but this
is immediate since $\det [x,y]_q = m_q \cdot m_q^{-1} \cdot m_q = m_q$.

Commutativity and associativity follow as in the proof of
\ref{algebra1}.  By assumption $\mathcal C(q)$ contains at least one
rational point.  Call this point $Q_0$ and denote $R_0 = Q_0$; clearly
$P \cdot Q_0^* \cdot R_0 = P$ for all points $P \in \mathcal C$.  The
equation 
$P \cdot Q_0^* \cdot R = R_0$ is satisfied by $R = Q_0 \cdot
P^* \cdot R_0$, proving that inverses do indeed exist. \end{proof}

This proposition motivates a definition of a mulltiplicative structure
similar to that of ternary algebras.

\begin{ternary2} A set $\mathcal C$ is called a commutative ternary group
if for $P,Q,R,S,T \in \mathcal C$, there is a triple product $P \cdot Q^*
\cdot R$ such that
\begin{enumerate}
   \item Multiplication: $P \cdot Q^* \cdot R \in \mathcal C$
   \item Commutativity: $P \cdot Q^* \cdot R = R \cdot Q^* \cdot P$
   \item Associativity: $\left( P \cdot Q^* \cdot R \right) \cdot
   S^* \cdot T = P \cdot Q^* \cdot \left( R \cdot S^* \cdot T \right)$
   \item Identity: There exist $Q_0, R_0 \in \mathcal C$ such that $P
   \cdot Q_0^* \cdot R_0 = P$ for all $P$.
   \item Inverses: Given $P$, there is $R$ such that $P \cdot Q_0^* \cdot R
= R_0$. 
\end{enumerate}
\end{ternary2}

\begin{fractions2} Let $\mathcal C$ be a commutative ternary group.  Denote
$\mathcal 
G$ as the collection of formal symbols $P \cdot Q^*$ for $P,Q \in \mathcal
C$ modulo the equivalence relation $P \cdot Q^* \sim R \cdot S^*$
whenever $P \cdot Q^* \cdot R_0 = R \cdot S^* \cdot R_0$ in $\mathcal
C$.

Then $\mathcal G$ is an abelian group with a transitive left
action on $\mathcal C$. \end{fractions2}

\begin{proof} The proof follows exactly as for \ref{fractions1}, so the
only statement we must prove is that for the group action.  Let $P, Q \in
\mathcal C$ be arbitrary; we must show that there is $g \in \mathcal G$
such that $P = g \cdot Q$.  Choose $R$ so that $Q \cdot Q_0^* \cdot R =
R_0$, and let $g = (P \cdot Q_0^* \cdot R) \cdot Q_0^*$.  Then
\[ g \cdot Q = (P \cdot Q_0^* \cdot R) \cdot Q_0^* \cdot Q = P \cdot
Q_0^* \cdot (R \cdot Q_0^* \cdot Q) = P \cdot Q_0^* \cdot R_0 = P \]

\noindent which completes the proof. \end{proof}


\section{Multiplication of Quadratic Forms}

The determinant maps $\mathcal A(q), \, \mathcal R(q) \to \mathbb Q$ now
give simple multiplicative formulas.

\begin{cor1} \label{cor1} Let $a$, $b$, and $c$ be integers.  The product
\[ \left[ a \, x_1^2 + b \, x_1 \, y_1 + c \, y_1^2 \right] \left[ a \,
x_2^2 + b \, x_2 \, y_2 + c \, y_2^2 \right] = u^2 + b \, u \, v + a \, c \,
v^2 \]

\noindent where $u = a \, x_1 \, x_2 + b \, x_1 \, y_2 + c \, y_1 \,
y_2$ and $v = y_1 \, x_2 - x_1 \, y_2$.  The triple product
\[ \left[ a \, x_1^2 + b \, x_1 \, y_1 + c \, y_1^2 \right] \left[ a \,
x_2^2 + b \, x_2 \, y_2 + c \, y_2^2 \right] \left[ a \, x_3^2 + b \,
x_3 \, y_3 + c \, y_3^2 \right] \]

\noindent is also in the form $a \, x^2 + b \, x \, y + c \, y^2$ where
\[ \begin{aligned}
   x & = a \, (x_1 \, x_2 \, x_3) + b \, (x_1 \, y_2 \, x_3) + c \,
   (x_1 \, y_2 \, y_3 - y_1 \, x_2 \, y_3 + y_1 \, y_2 \, x_3) \\
   y & = a \, (x_1 \, x_2 \, y_3 - x_1 \, y_2 \, x_3 + y_1 \, x_2 \,
   x_3) + b \, (y_1 \, x_2 \, y_3) + c \, (y_1 \, y_2 \, y_3)
\end{aligned} \]
\end{cor1}

This result is a generalization of the multiplicative formulas found
in \cite{MR90m:11016}.

\begin{proof} Set $q(x,y) = a \, x^2 + b \, x \, y + c \, y^2$, and first
assume that 
$\text{Disc}(q) \neq 0$.  Then $\text{Det}(q) = h_q = k_q = m_q = 0$,
$\mathcal 
A(q) = \left \{ [x,y]_q = \left[ \begin{matrix} x & -c \, y \\ y & a \,
x + b \, y \end{matrix} \right] \right \}$ and $\mathcal R(q) = \left \{
[u,v]_{q,0} = \left[ \begin{matrix} u & -c \, v \\ a \, v & u + b \, v
\end{matrix} \right] \right \}$.  Choose $A=[x_1,y_1]_q$,
$B=[x_2,y_2]_q$, and $C=[x_3,y_3]_q$ in $\mathcal A(q)$.  By \ref{algebra1},
$[u,v]_{q,0} = A 
\cdot B^* \in \mathcal R(q)$ and $[x,y]_q = A \cdot B^*
\cdot C \in \mathcal A(q)$.  Upon taking determinants of both sides we find
$u^2 + b \, 
u \, v + a \, c \, v^2 = \det A \cdot \det B$ and $a \, x^2 +
b \, x \, y + c \, y^2 = \det A \cdot \det B \cdot \det C$.

The formulas for $u$, $v$, $x$ and $y$ hold even for $\text{Disc}(q) =
0$ so they are always valid.  \end{proof}

As another application, we state without proof the following.

\begin{cor2} Let $(x_1,y_1,z_1)$, $(x_2,y_2,z_2)$, and $(x_3,y_3,z_3)$
be points on the projective variety $a \, x^2 + b \, y^2 + c \, z^2 =
0$.  A fourth point is $(x,y,z)$ given by
\[ \begin{aligned}
   x & = a \, (x_1 \, x_2 \, x_3) + b \, (x_1 \, y_2 \, y_3 - y_1 \,
   x_2 \, y_3 + y_1 \, y_2 \, x_3)  \\
   y & = a \, (x_1 \, x_2 \, y_3 - x_1 \, y_2 \, x_3 + y_1 \, x_2 \,
   x_3) + b \, (y_1 \, y_2 \, y_3) \\
   z & = c \, (z_1 \, z_2 \, z_3)
\end{aligned} \]
\end{cor2}

We conclude with an application to the values of a given binary
quadratic form.

\begin{algebra3} Fix $q(x,y)$ a binary quadratic form, define $\mathbb
Q(q)$ as the set of $\alpha = q(x,y)$ for $x,y \in \mathbb Q$ such that
$\text{Disc}(q) \cdot \alpha \neq \text{Det}(q)$, and assume that $\mathbb
Q(q)$ is nonempty.
Then for $\alpha, \beta, \gamma \in \mathbb Q(q)$, the map
\[ \alpha \cdot \beta^* \cdot \gamma = \frac {\text{Disc}(q) \, \alpha
\, \gamma - \text{Det}(q) \, (\alpha - \beta + \gamma)}{\text{Disc}(q)
\, \beta - \text{Det}(q)} \]

\noindent turns $\mathbb Q(q)$ into a commutative ternary group.

\end{algebra3}

$\mathbb Q(q)$ is the image of $q(x,y)$.  When $\text{Det}(q) = 0$, such
as with $q(x,y) = x^2 + y^2$, we recover formulas similar to those in
\ref{cor1} because the triple product is $\alpha \cdot \beta^{-1} \cdot
\gamma$.  When $\text{Disc}(q) = 0$, such as with $q(x,y) = x^2 - y$, the
``product'' is $\alpha - \beta + \gamma$.  Projectively the parabola
``looks'' like the line so it makes sense that ``multiplication'' in this
case is really addition.

\begin{proof}  First assume that $\text{Disc}(q) \neq 0$.  Choose $\alpha =
q(P)$, $\beta=q(Q)$, and $\gamma=q(R)$
in $\mathbb Q(q)$ for some points $P=(x_1,y_1)$, $Q=(x_2,y_2)$ and
$R=(x_3,y_3)$.  We map them to matrices $A=[x_1,y_1]_q$,
$B=[x_2,y_2]_q$, and $C=[x_2,y_3]_q$ in $\mathcal A(q)$ where $\det
A = \alpha + m_q$, $\det B = \beta + m_q$, and $\det C = \gamma + m_q$
are all nonzero by assumption since $m_q = -\text{Det}(q)/\text{Disc}(q)$.
Then the triple product $[x,y]_q = A \cdot B^{-1} \cdot C$ exists so we
define
\[ \alpha \cdot \beta^* \cdot \gamma = q(x,y) = \det [x,y]_q - m_q = \frac
{(\alpha + m_q)(\gamma + m_q)}{\beta + m_q} - m_q \]

\noindent It follows from a proof similar to that in \ref{algebra2} that
$\mathbb Q(q)$ is a commutative ternary group.

Now we consider the case where $\text{Disc}(q) = 0$.  Since $\mathbb
Q(q)$ is nonempty by assumption we have $\text{Det}(q) \neq 0$.  The
formulas 
for $x$ and $y$ reduce to $(x,y) = P - Q + R$ --- even though the ternary
algebra $\mathcal A(q)$ is not defined --- and the ternary group operation
reduces to $\alpha \cdot \beta^* \cdot \gamma = \alpha - \beta + \gamma$.
\end{proof}


\bibliographystyle{plain}

\end{document}